\documentclass[11pt,preprint]{article}
\usepackage{mathrsfs,amsmath,amssymb,epsfig,subfigure,color}

\renewcommand{\arraystretch}{1.2}

\newcommand{\BE}{\begin{equation}}
\newcommand{\EE}{\end{equation}}

\newtheorem{theorem}{Theorem}[section]
\newtheorem{lemma}{Lemma}[section]
\newtheorem{example}{Example}[section]

\newtheorem{remark}{Remark} [section]

\def\dfrac{\displaystyle\frac}
\def\dsum{\displaystyle\sum}
\def\dint{\displaystyle\int}

\begin{document}
\title{\bf A high order ADI scheme for the two-dimensional time fractional diffusion-wave equation}
\author{Zhibo Wang\thanks{Corresponding author. Email: zhibowangok@gmail.com. Department of Mathematics, University of Macau, Av. Padre Tom\'{a}s Pereira Taipa, Macau, China.} \and Seakweng Vong\thanks{Email: swvong@umac.mo. Department of Mathematics, University of Macau, Av. Padre Tom\'{a}s Pereira Taipa, Macau, China.}}
\date{}
 \maketitle

 \begin{abstract}
 In this paper, a compact alternating direction implicit (ADI) finite
 difference scheme for the two-dimensional time fractional diffusion-wave
 equation is developed, with temporal and spatial accuracy order equal to
 two and four respectively. The second order accuracy in the time direction
 has not been achieved in previous studies.
 \end{abstract}

\section{Introduction}
 In this paper, we concentrate on high order compact ADI finite difference
 methods for the following two-dimensional fractional diffusion-wave (2DFDW) equation:
 \begin{align}\label{main1}
 &_{0}^CD_t^{\gamma}u
 =\Delta u+g(x,y,t), \quad (x,y)\in \Omega, \quad 0<t\leq T, \quad 1<\gamma<2,\\
 &u(x,y,0)=\psi(x,y), \quad \frac{\partial u(x,y,0)}{\partial t}=\phi(x,y), \quad (x,y)\in\bar{\Omega}=\Omega\cup\partial\Omega,\label{main2}\\
 &u(x,y,t)=\varphi(x,y,t), \quad (x,y)\in \partial\Omega, \quad 0<t\leq T,\nonumber
 \end{align}
 where $_{0}^CD_t^{\gamma}u$ is the Caputo fractional derivative of $u$
 with respect to the time variable $t$ defined as
 $$_{0}^CD_t^{\gamma}u(x,y,t)=\frac1{\Gamma(2-\gamma)}\int_0^t
 \frac{\partial^2u(x,y,s)}{\partial s^2}(t-s)^{1-\gamma}ds,$$
 with $\Gamma(\cdot)$ being the gamma function, $\Delta$ is the
 two-dimensional Laplacian, $\Omega=(0,L_1)\times(0,L_2)$ and $\partial\Omega$ is
 the boundary of $\Omega$. For detailed applications of this problem, we  refer readers to \cite{Mainardi1,Mainardi2}.

 Progress on the study of high order scheme for time fractional differential
 equations is rapid in recent years. We refer the readers to
 \cite{Gao,Cui,Mohebbi,Wang} and the references therein for details on the
 development of one-dimensional problems. The ADI method is efficient
 for solving high dimensional problems. This method has been successfully applied to time fractional differential equations
 \cite{TandM,Zhang2,Cuicompact1} in the past years. To our
 knowledge, the most recent result on ADI scheme for the 2DFDW equation
 with Dirichlet boundary condition is given by \cite{Zhang2}, where a
 compact scheme of order $(\tau^{3-\gamma}+h_1^4+h_2^4)$ was established.
 Here, $\tau$ is the temporal grid size and $h_1$, $h_2$ are the spatial
 grid sizes in different directions.

 In this paper, we introduce an ADI scheme for the 2DFDW equation. The
 scheme is shown to converge with order $(\tau^2+h_1^4+h_2^4)$. This
 scheme is based on our recent result \cite{Wang} on one-dimensional
 time fractional differential equations. The idea, inspired by
 \cite{Meerschaert,Tian} on the study of space fractional differential
 equations, is to carry out discretization using the shifted Gr\"{u}nwald
 difference operator. The details are given in the next two sections which are
 justified by numerical examples in Section 4.

 \section{The proposed compact ADI scheme}
 We assume that $\psi\equiv0$ in (\ref{main2}) without loss of generality
 since we can solve the equation for $v(x,y,t)=u(x,y,t)-\psi(x,y)$ in
 general.

 In previous studies, discretization is directly based on (\ref{main1}).
 However, we find that the accuracy of the discrete approximations can be improved
 if we consider an equivalent form of (\ref{main1}). For details of this equivalent form,
 we refer readers to \cite{Huang}. We just outline the key elements here for the completeness of our presentation.
  We first note that
 \begin{align*}
 {_{0}^CD_t^{\gamma}}u(x,y,t)&=\dfrac1{\Gamma(1-(\gamma-1))}
 \dint_0^t(t-s)^{-(\gamma-1)}\dfrac{\partial}{\partial s}\dfrac{\partial u(x,y,s)}{\partial s}ds\\
 &={_{0}^CD_t^{\gamma-1}}\dfrac{\partial u(x,y,t)}{\partial t}
 ={_{0}^CD_t^{\alpha}}\dfrac{\partial u(x,y,t)}{\partial t},
 \end{align*}
 where $0<\alpha=\gamma-1<1$. Taking the  Riemann-Liouville fractional integral $_{0}I^\alpha_t$ of order $\alpha$ on both sides of (\ref{main1}),
 we obtain
 \begin{equation}\label{changed-main}
 \dfrac{\partial u(x,y,t)}{\partial t}=\phi(x,y)+
 \dfrac1{\Gamma(\alpha)}\dint_0^t(t-s)^{\alpha-1}\Delta u(x,y,s)ds+f(x,y,t), \quad (x,y)\in\Omega, \quad 0<t\leq T,
 \end{equation}
 where $f(x,y,t)={_{0}}I^\alpha_tg(x,y,t)$ and
 the Riemann-Louville frational integral of a function $h(t)$ is defined by
 $$_{0}I^\alpha_th(t)=\dfrac1{\Gamma(\alpha)}\dint_0^t(t-s)^{\alpha-1}h(s)ds.$$

 We develop our scheme based on (\ref{changed-main}) by the weighted and shifted Gr\"{u}nwald difference formula for the Riemann-Liouville fractional integral. To this end, we let $h_1=\frac{L_1}{M_1},~h_2=\frac{L_2}{M_2}$ and $\tau=\frac TN$ be the spatial and temporal step sizes respectively, where  $M_1,~M_2$ and $N$ are some given integers.
 For $i=0,1,\ldots,M_1,~j=0,1,\ldots,M_2$ and $k=0,1,\ldots,N$, denote $x_i=ih_1,~y_j=jh_2,~t_k=k\tau$.
  In this paper, we consider the problem under the Dirchlet boundary conditions and
  we need to determine the approximated values of the solution on
  the $(M_1-1)\times(M_2-1)$ interior grid points. For any grid function $u=\{u_{ij}^k|0\leq i\leq M_1,~0\leq j\leq M_2,~0\leq k\leq N\}$, we introduce the following notations:
 $$\begin{array}{c}
 u_{ij}^{k-\frac12}=\frac12(u_{ij}^{k}+u_{ij}^{k-1}),
 ~~\delta_xu_{i-\frac12,j}^k=\frac1{h_1}(u_{ij}^k-u_{i-1,j}^k),
 ~~\delta_x^2u_{ij}^k=\frac1{h_1}(\delta_xu_{i+\frac12,j}^k-\delta_xu_{i-\frac12,j}^k),\\[3pt]
 {\cal H}_xu_{ij}=\left\{
 \begin{array}{ll}
 \Big(1+\frac{h_1^2}{12}\delta_x^2\Big)u_{ij}
 =\frac1{12}(u_{i-1,j}+10u_{ij}+u_{i+1,j}),
 &1\leq i\leq M_1-1,~0\leq j\leq M_2,\\
 u_{ij}, &i=0 \mbox{ or } M_1,~0\leq j\leq M_2.
 \end{array}\right.
 \end{array}$$
 Corresponding notations in the $y$ direction are defined similarly.
 We further denote
 $$\begin{array}{c}
 \delta_y\delta_xu_{i-\frac12,j-\frac12}^k
 =\frac1{h_2}(\delta_xu_{i-\frac12,j}^k-\delta_xu_{i-\frac12,j-1}^k),
 ~~\delta_y\delta_x^2u_{i,j-\frac12}^k
 =\frac1{h_2}(\delta_x^2u_{ij}^k-\delta_x^2u_{i,j-1}^k),\\[4pt]
 {\cal H}u_{ij}={\cal H}_x{\cal H}_yu_{ij},~~\Lambda u_{ij}=({\cal H}_y\delta^2_x+{\cal H}_x\delta_y^2)u_{ij},~~\langle u,v\rangle=h_1h_2\sum\limits_{i=1}^{M_1-1}\sum\limits_{j=1}^{M_2-1}u_{ij}v_{ij},\\
 \langle \delta_xu,\delta_xv\rangle
 =h_1h_2\sum\limits_{i=1}^{M_1}\sum\limits_{j=1}^{M_2-1}\delta_xu_{i-\frac12,j}\delta_xv_{i-\frac12,j},
 ~~\langle\delta_x\delta_yu,\delta_x\delta_yv\rangle=h_1h_2\sum\limits_{i=1}^{M_1}\sum\limits_{j=1}^{M_2}
 \delta_x\delta_yu_{i-\frac12,j-\frac12}\delta_x\delta_yv_{i-\frac12,j-\frac12},\\
 \|u\|^2=\langle u,u\rangle,~~\|\delta_xu\|^2=\langle \delta_xu,\delta_xv\rangle,
 ~~\|\delta_x\delta_yu\|^2=\langle\delta_x\delta_yu,\delta_x\delta_yu\rangle.
 \end{array}$$

 The high order accuracy of our proposed scheme is based on the second order approximation of the
 Riemann-Liouville fractional integral. This approximation is established by
 the shifted operator for the Riemann-Liouville fractional integral defined as:
 $$\begin{array}{c}
 \mathcal{A}_{\tau,r}^\alpha f(t)=\tau^\alpha\sum\limits_{k=0}^{\infty}\omega_kf(t-(k-r)\tau),
 \end{array}$$
 where $\omega_k=(-1)^k\binom{-\alpha}{k}$.
 By using Fourier transform and ideas in \cite{Tian}, we obtained the following estimate  in \cite{Wang}.
 \begin{lemma}\label{main-operator}
 {Assume that $f(t),~{_{-\infty}I_t^\alpha f(t)}$ and
 $(i\omega)^{2-\alpha}\mathscr{F}[f](\omega)$} belong to $L^1(\mathbb{R})$.
 Define the weighted and shifted difference operator by
 $$\mathcal{I}_{\tau,p,q}^\alpha f(t)
 =\dfrac{2q+\alpha}{2(q-p)}\mathcal{A}_{\tau,p}^\alpha f(t)
 +\dfrac{2p+\alpha}{2(p-q)}\mathcal{A}_{\tau,q}^\alpha f(t).$$
 Then we have
 $$\mathcal{I}_{\tau,p,q}^\alpha f(t)={_{-\infty}I_t^\alpha f(t)}+O(\tau^2)$$
 for $t\in \mathbb{R}$, where $p$ and $q$ are integers and $p\neq q$.
 \end{lemma}
 {\bf Proof.} The proof can be found in \cite{Wang}.
 However, since this lemma plays a crucial role for our main result,
 we restate a simplified version here for the completeness of our presentation.

 Referring to the definition of $\mathcal{A}_{\tau,r}^\alpha$, we let
 \begin{equation}\label{weight}
 \mathcal{I}_{\tau,p,q}^\alpha f(t)=\tau^\alpha\bigg[\mu_1\dsum_{k=0}^{\infty}\omega_k^{(\alpha)}f(t-(k-p)\tau)
 +\mu_2\dsum_{k=0}^{\infty}\omega_k^{(\alpha)}f(t-(k-q)\tau)\bigg].
 \end{equation}
 The basic idea of obtaining high order accuracy is to apply Fourier transform to (\ref{weight}), which yields
 \begin{equation}\label{fourior-weight}
 \mathscr{F}[\mathcal{I}_{\tau,p,q}^\alpha f](\omega)
 =(i\omega)^{-\alpha}[\mu_1W_p(i\omega\tau)+\mu_2W_q(i\omega\tau)]\mathscr{F}[f](\omega),
 \end{equation}
 where
 \begin{equation}\label{taylor}
 W_r(z)=z^\alpha(1-e^{-z})^{-\alpha}e^{rz}=1+\big(r+\frac\alpha2\big)z+O(z^2),
 \quad r=p,~q.
 \end{equation}
 Notice that \cite{Ervin}
 $$\mathscr{F}[{_{-\infty}I_t^\alpha f(t)}]=(i\omega)^{-\alpha}\mathscr{F}[f](\omega).$$
 From (\ref{fourior-weight}) and (\ref{taylor}), one can see that,
 by taking $\mu_1=\frac{2q+\alpha}{2(q-p)}$ and
 $\mu_2=\frac{2p+\alpha}{2(p-q)}$,
 the constant and linear term in $\mathscr{F}[\mathcal{I}_{\tau,p,q}^\alpha f-{_{-\infty}I_t^\alpha f}]$ can be eliminated.
 We can therefore conclude that
 $$|\mathcal{I}_{\tau,p,q}^\alpha f-{_{-\infty}I_t^\alpha f}|
 \leq\dfrac1{2\pi}\dint_\mathbb{R}|\mathscr{F}[\mathcal{I}_{\tau,p,q}^\alpha f-{_{-\infty}I_t^\alpha f}]|d\omega
 \leq C\|(i\omega)^{2-\alpha}\mathscr{F}[f](\omega)\|_{L^1}
 \tau^2
 =O(\tau^2).\qquad\Box$$

 Since the solution $u(x,y,t)$  can be continuously extended to be zero for $t<0$,
 by choosing $(p,q)=(0,-1)$, which yields $\frac{2q+\alpha}{2(q-p)}=1-\frac\alpha2,
 ~\frac{2p+\alpha}{2(p-q)}=\frac\alpha2$, one can obtain
 \begin{equation}\label{rl-lap}
 \begin{array}{rl}
 _{0}I^\alpha_t\Delta u(x_i,y_j,t_{n+1})
 &=\tau^\alpha\bigg[\Big(1-\dfrac\alpha2\Big)\dsum_{k=0}^{n+1}\omega_k(\delta^2_x+\delta^2_y)u^{n+1-k}_{ij}
 +\dfrac\alpha2\dsum_{k=0}^{n}\omega_k(\delta^2_x+\delta^2_y)u^{n-k}_{ij}\bigg]
 +(R_1)_{ij}^{n+1}\\
 &=\tau^\alpha\dsum_{k=0}^{n+1}\lambda_k(\delta^2_x+\delta^2_y)u^{n+1-k}_{ij}
 +(R_1)_{ij}^{n+1},
 \end{array}
 \end{equation}
 where $(R_1)_{ij}^{n+1}=O(\tau^2+h_1^2+h_2^2)$ and
 \begin{equation}\label{lambda}
 \lambda_0=(1-\dfrac\alpha2)\omega_0,
 ~\lambda_k=(1-\dfrac\alpha2)\omega_{k}+\dfrac\alpha2\omega_{k-1},~k\geq1.
 \end{equation}

 To raise the accuracy in spatial directions, we need
 \begin{lemma}\label{compact}{\rm(\cite{Gao})}
 If $f(x)\in {\cal C}^6[x_{i-1},x_{i+1}],~1\leq i\leq M_1-1$, then it holds that
 $$\frac1{12}[f''(x_{i-1})+10f''(x_i)+f''(x_{i+1})]
 =\frac1{h_1^2}[f(x_{i-1})-2f(x_i)+f(x_{i+1})]+O(h_1^4).$$
 \end{lemma}

 Denoting $\mu=\frac{\tau^{\alpha+1}}2$, Lemma
 \ref{compact} and (\ref{rl-lap}) yield the following weighted compact Crank-Nicolson scheme:
 \begin{equation}\label{compact1}
 {\cal H}(u^{n+1}_{ij}-u^n_{ij})=\tau{\cal H}\phi_{ij}
 +\mu\bigg(\dsum_{k=0}^{n+1}\lambda_k\Lambda u^{n+1-k}_{ij}
 +\dsum_{k=0}^{n}\lambda_k\Lambda u^{n-k}_{ij}\bigg)
 +\dfrac\tau2{\cal H}(f^n_{ij}+f^{n+1}_{ij})+\tau(R_2)_{ij}^{n+1},
 \end{equation}
 where $(R_2)_{ij}^{n+1}=O(\tau^2+h_1^4+h_2^4)$.

 \medskip
 Adding a small term
 $\mu^2\lambda_0^2\delta_x^2\delta_y^2(u_{ij}^{n+1}-u_{ij}^{n})=O(\tau^{3+2\alpha})$
 on both sides of (\ref{compact1}), we have
 \begin{equation}\label{compact-scheme}
 \begin{array}{rl}
 &\quad{\cal H}(u^{n+1}_{ij}-u^n_{ij})+
 \mu^2\lambda_0^2\delta_x^2\delta_y^2(u_{ij}^{n+1}-u_{ij}^{n})\\
 &=\tau{\cal H}\phi_{ij}
 +\mu\bigg(\dsum_{k=0}^{n+1}\lambda_k\Lambda u^{n+1-k}_{ij}
 +\dsum_{k=0}^{n}\lambda_k\Lambda u^{n-k}_{ij}\bigg)
 +\dfrac\tau2{\cal H}(f^n_{ij}+f^{n+1}_{ij})+\tau R_{ij}^{n+1}
 \end{array}
 \end{equation}
 with $R_{ij}^{k+1}=O(\tau^{2}+h_1^4+h_2^4)$. Omitting the truncation error in (\ref{compact-scheme}), we reach the following ADI scheme
 $$\begin{array}{ll}
 \quad({\cal H}_x-\mu\lambda_0\delta^2_x)
 ({\cal H}_y-\mu\lambda_0\delta_y^2)u_{ij}^{n+1}\\
 =({\cal H}_x+\mu\lambda_0\delta^2_x)
 ({\cal H}_y+\mu\lambda_0\delta_y^2)u_{ij}^{n}
 +\mu\bigg(\dsum_{k=1}^{n+1}\lambda_k\Lambda u^{n+1-k}_{ij}
 +\dsum_{k=1}^{n}\lambda_k\Lambda u^{n-k}_{ij}\bigg)
 +\tau{\cal H}\phi_{ij}+\dfrac\tau2{\cal H}(f^n_{ij}+f^{n+1}_{ij}),\\
 \quad(x_i,y_j)\in\Omega, \quad 1\leq n\leq N-1,\\
 u_{ij}^0=0, \quad (x_i,y_j)\in\bar{\Omega},\\
 u_{ij}^n=\varphi(x_i,y_j,t_n), \quad 1\leq n\leq N, \quad (x_i,y_j)\in\partial\Omega.
 \end{array}$$

 We note that, in ADI methods (see \cite{Zhang2} for example),
 the solution $\{u_{ij}^{n+1}\}$ is determined by solving two independent one-dimensional problems. Namely, the intermediate variables
 $$u_{ij}^*=({\cal H}_y-\mu\lambda_0\delta_y^2)u_{ij}^{n+1},~~1\leq i\leq M_1-1,~1\leq j\leq M_2-1,$$
 are first solved from the following system with fixed $j\in\{1,2,\ldots,M_2-1\}$:
 $$\left\{\begin{array}{ll}
 ({\cal H}_x-\mu\lambda_0\delta^2_x)u_{ij}^*
 =({\cal H}_x+\mu\lambda_0\delta^2_x)
 ({\cal H}_y+\mu\lambda_0\delta_y^2)u_{ij}^{n}
 +\mu\bigg(\dsum_{k=1}^{n+1}\lambda_k\Lambda u^{n+1-k}_{ij}
 +\dsum_{k=1}^{n}\lambda_k\Lambda u^{n-k}_{ij}\bigg)\\
 \qquad\qquad\qquad\qquad~+\tau{\cal H}\phi_{ij}
 +\dfrac\tau2{\cal H}(f^n_{ij}+f^{n+1}_{ij}),~~1\leq i\leq M_1-1,\\[3pt]
 u_{0j}^*=({\cal H}_y-\mu\lambda_0\delta_y^2)u_{0j}^{n+1}
 =({\cal H}_y-\mu\lambda_0\delta_y^2)[\varphi(x_0,y_j,t_{n+1})],\\
 u_{M_1j}^*=({\cal H}_y-\mu\lambda_0\delta_y^2)u_{M_1j}^{n+1}
 =({\cal H}_y-\mu\lambda_0\delta_y^2)[\varphi(x_{M_1},y_j,t_{n+1})].
 \end{array}\right.$$
 Once $\{u_{ij}^*\}$ is available, we then solve $\{u_{ij}^{n+1}\}$ from the following system for fixed $i\in\{1,2,\ldots,M_1-1\}$:
 $$\left\{\begin{array}{ll}
 ({\cal H}_y-\mu\lambda_0\delta_y^2)u_{ij}^{n+1}=u_{ij}^*,
 ~~1\leq j\leq M_2-1,\\[3pt]
 u_{i0}^{n+1}=\varphi(x_i,y_0,t_{n+1}),
 ~~u_{iM_2}^{n+1}=\varphi(x_i,y_{M_2},t_{n+1}).
 \end{array}\right.$$

 In the next section, the convergence result will be proved.

 \section{Convergence analysis of the compact ADI scheme}
 We first introduce two useful lemmas.

 \medskip
 \begin{lemma}\label{sequence}{\rm(\cite{Wang})}
 Let $\{\lambda_n\}_{n=0}^\infty$ be defined as {\rm(\ref{lambda})}, then for any positive integer $k$ and
 $(v_1,v_2,\ldots,v_k)^T\in\mathbb{R}^k$,
 it holds that
 $$\dsum_{n=0}^{k-1}\Big(\dsum_{p=0}^n\lambda_pv_{n+1-p}\Big)v_{n+1}\geq0.$$
 \end{lemma}
 \begin{lemma}\label{gronwall}{\rm(Grownall's inequality \cite{Quarteroni})}
 Assume that $\{k_n\}$ and $\{p_n\}$ are nonnegative sequences, and the sequence $\{\phi_n\}$ satisfies
 $$\phi_0\leq g_0,~~~\phi_n\leq g_0+\dsum_{l=0}^{n-1}p_l
 +\dsum_{l=0}^{n-1}k_l\phi_l,~~~n\geq1,$$
 where $g_0\geq0$. Then the sequence $\{\phi_n\}$ satisfies
 $$\phi_n\leq\Big(g_0+\dsum_{l=0}^{n-1}p_l\Big)
 \exp{\Big(\dsum_{l=0}^{n-1}k_l\Big)},~~~n\geq1.$$
 \end{lemma}

 With the above lemmas, we now proceed to prove the convergence of our compact ADI scheme (\ref{compact-scheme}).
 \begin{theorem}\label{convergence-th}
 Assume that $u(x,y,t)\in {\cal C}_{x,y,t}^{6,6,2}(\Omega\times[0,T])$ is the solution of {\rm(\ref{changed-main})} and $\{u_{ij}^k|0\leq i\leq M_1,~ 0\leq j\leq M_2,~ 0\leq k\leq N\}$
 is a solution of the finite difference scheme {\rm(\ref{compact-scheme})}, respectively. Denote
 $$e_{ij}^k=u(x_i,y_j,t_k)-u_{ij}^k,~~(x_i,y_j)\in\bar{\Omega},~~0\leq k\leq N.$$
 Then there exists a positive constant $\tilde c$ such that
 $$\|e^k\|\leq\tilde c(\tau^2+h_1^4+h_2^4), \quad 0\leq k\leq N.$$
 \end{theorem}
{\bf Proof.} We can easily get the following error equation
 \begin{align}\label{error}
 &{\cal H}(e_{ij}^{k+1}-e_{ij}^{k})
 +\mu^2\lambda_0^2\delta_x^2\delta_y^2(e_{ij}^{k+1}-e_{ij}^{k})
 =\mu\dsum_{l=0}^{k}\lambda_l\Lambda(e_{ij}^{k+1-l}+e_{ij}^{k-l})
 +\tau R_{ij}^{k+1},\\
 &\quad1\leq i\leq M_1-1,~1\leq j\leq M_2-1,~0\leq k\leq N-1,\nonumber\\
 &e_{ij}^0=0, \quad (x_i,y_j)\in\bar{\Omega},\nonumber\\
 &e_{ij}^k=0, \quad 1\leq k\leq N, \quad (x_i,y_j)\in\partial\Omega.\nonumber
 \end{align}
 Multiplying (\ref{error}) by
 $h_1h_2(e_{ij}^{k+1}+e_{ij}^{k})$ and summing in $i,j$, we obtain
 \begin{equation}\label{error-inner}
 \begin{array}{rl}
 &\big\langle{\cal H}(e^{k+1}-e^{k}),(e^{k+1}+e^{k})\big\rangle
 +\mu^2\lambda_0^2\big\langle\delta_x^2\delta_y^2(e^{k+1}-e^{k})
 ,(e^{k+1}+e^{k})\big\rangle\\
 =&4\mu\dsum_{l=0}^{k}\lambda_l\big\langle\Lambda e^{k+\frac12-l},e^{k+\frac12}\big\rangle
 +\tau\big\langle R^{k+1},(e^{k+1}+e^{k})\big\rangle,\\
 &1\leq i\leq M_1-1,~1\leq j\leq M_2-1,~0\leq k\leq N-1.
 \end{array}
 \end{equation}
 Note that the term $\big\langle{\cal H}(e^{k+1}-e^{k}),(e^{k+1}+e^{k})\big\rangle$ on the left
 can be written as $\big\langle{\cal H}e^{k+1},e^{k+1}\big\rangle-\big\langle{\cal H}e^{k},e^{k}\big\rangle$.

 Notice that, if $v$, $w$ satisfy   $v_{0,j}=v_{M_1,j}=0$ for all $j$ and
 $w_{i,0}=w_{i,M_2}=0$ for all $i$, we have
 $\langle\delta_x^2u,v\rangle=-\langle\delta_xu,\delta_xv\rangle$
 and $\langle\delta_y^2u,w\rangle=-\langle\delta_yu,\delta_yw\rangle$.
 The second term in the left hand side can therefore  be written as
 $$\mu^2\lambda_0^2\big\langle\delta_x^2\delta_y^2(e^{k+1}-e^{k})
 ,(e^{k+1}+e^{k})\big\rangle
 =\mu^2\lambda_0^2\Big(\big\langle\delta_x^2\delta_y^2e^{k+1},
 e^{k+1}\big\rangle
 -\big\langle\delta_x^2\delta_y^2e^{k},e^{k}\big\rangle\Big).$$
 Now we turn to the first term in the right hand side.
 Since ${\cal H}_x$ and ${\cal H}_y$ are positive definite and
 self-adjoint, we can consider their square roots denoted as
 ${\cal Q}_x$ and ${\cal Q}_y$ respectively.
 We thus have
 $$\begin{array}{rl}
 \big\langle\Lambda e^{k+\frac12-l},e^{k+\frac12}\big\rangle
 &=\big\langle{\cal H}_y\delta^2_xe^{k+\frac12-l},e^{k+\frac12}\big\rangle
 +\big\langle{\cal H}_x\delta^2_ye^{k+\frac12-l},e^{k+\frac12}\big\rangle\\[3pt]
 &=-\big\langle{\cal H}_y\delta_xe^{k+\frac12-l},
 \delta_xe^{k+\frac12}\big\rangle
 -\big\langle{\cal H}_x\delta_ye^{k+\frac12-l},
 \delta_ye^{k+\frac12}\big\rangle\\[3pt]
 &=-\big\langle {\cal Q}_y\delta_xe^{k+\frac12-l},
 {\cal Q}_y\delta_xe^{k+\frac12}\big\rangle
 -\big\langle {\cal Q}_x\delta_ye^{k+\frac12-l},
 {\cal Q}_x\delta_ye^{k+\frac12}\big\rangle.
 \end{array}$$
 Summing up (\ref{error-inner}) for $0\leq k\leq n-1$, we have
 $$\begin{array}{ll}
 \quad\langle{\cal H}e^{n},e^n\rangle
 +\mu^2\lambda_0^2\big\langle\delta_x^2\delta_y^2e^{n},e^{n}\big\rangle\\
 \leq-4\mu\dsum_{k=0}^{n-1}\dsum_{l=0}^{k}\lambda_l
 \Big(\big\langle {\cal Q}_y\delta_xe^{k+\frac12-l},
 {\cal Q}_y\delta_xe^{k+\frac12}\big\rangle
 +\big\langle {\cal Q}_x\delta_ye^{k+\frac12-l},
 {\cal Q}_x\delta_ye^{k+\frac12}\big\rangle\Big)
 +\tau\dsum_{k=0}^{n-1}\langle R^{k+1},(e^{k+1}+e^{k})\rangle.
 \end{array}$$
  Notice $\big\langle\delta_x^2\delta_y^2e^{n},e^{n}\big\rangle
 =\|\delta_x\delta_ye^n\|^2\geq0$ and,
 by $\|\delta_xu\|^2\leq\frac4{h_1^2}\|u\|^2$, we deduce
 $$\begin{array}{rl}
 \langle{\cal H}e^{n},e^n\rangle
 &=\langle e^{n},e^n\rangle+\dfrac{h_1^2}{12}\langle\delta_x^2e^{n},e^n\rangle
 +\dfrac{h_2^2}{12}\langle\delta_y^2e^{n},e^n\rangle
 +\dfrac{h_1^2h_2^2}{144}\langle\delta_x^2\delta_y^2e^{n},e^n\rangle\\[5pt]
 &=\langle e^{n},e^n\rangle-\dfrac{h_1^2}{12}\langle\delta_xe^{n},\delta_xe^n\rangle
 -\dfrac{h_2^2}{12}\langle\delta_ye^{n},\delta_ye^n\rangle
 +\dfrac{h_1^2h_2^2}{144}\langle\delta_x\delta_ye^{n},\delta_x\delta_ye^n\rangle\\[5pt]
 &\geq\dfrac13\|e^n\|^2.
 \end{array}$$
 It then follows from Lemma \ref{sequence} that
 $$\begin{array}{rl}
 \dfrac13\|e^n\|^2
 &\leq\tau\dsum_{k=0}^{n-1}\langle R^{k+1},(e^{k+1}+e^{k})\rangle\\[1pt]
 &\leq\dfrac14\|e^n\|^2+\tau^2\|R^{n}\|^2+\dfrac\tau4\|e^{n-1}\|^2+\tau\|R^{n}\|^2
 +\dfrac\tau4\dsum_{k=1}^{n-1}\|e^{k}\|^2
 +\dfrac\tau4\dsum_{k=1}^{n-2}\|e^{k}\|^2
 +2\tau\dsum_{k=0}^{n-2}\|R^{k+1}\|^2\\
 &\leq\dfrac14\|e^n\|^2+\tau^2\|R^{n}\|^2
 +\dfrac{\tau}2\dsum_{k=1}^{n-1}\|e^{k}\|^2
 +2\tau\dsum_{k=0}^{n-1}\|R^{k+1}\|^2,
 \end{array}$$
 which gives
 \begin{equation}\label{for-stability}
 \|e^n\|^2\leq12\tau^2\|R^{n}\|^2
 +6\tau\dsum_{k=1}^{n-1}\|e^{k}\|^2
 +24\tau\dsum_{k=0}^{n-1}\|R^{k+1}\|^2
 \leq6\tau\dsum_{k=1}^{n-1}\|e^{k}\|^2+c(\tau^2+h_1^4+h_2^4)^2,
 \end{equation}
 the desired result then follows by Lemma \ref{gronwall}.$\qquad\Box$

 \begin{remark}
 Following the idea of the proof for Theorem \ref{convergence-th},
 one can  show that the proposed compact ADI scheme {\rm(\ref{compact-scheme})} is unconditionally stable.
 In fact, assume that
 $u _{ij}^k=0$ for $(x_i,y_j)\in\partial\Omega$
 and $\{u_{i,j}^n\}$ is the solution of the scheme
 $${\cal H}(u^{n+1}_{ij}-u^n_{ij})+
 \mu^2\lambda_0^2\delta_x^2\delta_y^2(u_{ij}^{n+1}-u_{ij}^{n})
 =\mu\bigg(\dsum_{k=0}^{n+1}\lambda_k\Lambda u^{n+1-k}_{ij}
 +\dsum_{k=0}^{n}\lambda_k\Lambda u^{n-k}_{ij}\bigg)
 +\tau{\cal H}q_{ij}+\tau g^{n+\frac12}_{ij},$$
 where $(\frac{\partial u}{\partial t})_{ij}^0=q_{ij},~1\leq i\leq M_1-1,~1\leq j\leq M_2-1,~k\geq0$.
 Then, by arguments similar to those for getting {\rm(\ref{for-stability})}, we can obtain
 \begin{align*}
 \|u^n\|^2&\leq12\|u^0\|^2+6\tau\dsum_{k=0}^{n-1}\|u^{k}\|^2
 +25\tau\dsum_{k=0}^{n-1}\|\Lambda u^0+{\cal H}q+g^{k+\frac12}\|^2\\
 &\leq e^{6T}\Big(12\|u^0\|^2+25\tau\dsum_{k=0}^{n-1}\|\Lambda u^0+{\cal H}q+g^{k+\frac12}\|^2\Big).
 \end{align*}
 \end{remark}

 \section{Numerical experiments}
 In this section, we carry out numerical experiments for our compact finite difference scheme.
 \begin{example}\label{ex1}
  The following problem is considered in {\rm\cite{Zhang2}}:
 $$\begin{array}{rl}
 &_{0}^CD_t^{\gamma}u=\Delta u
 +\sin(x)\sin(y)\big[\frac{\Gamma(\gamma+3)}2t^2+2t^{\gamma+2}\big],
 ~~(x,y)\in\Omega=(0,\pi)\times(0,\pi),~~0<t\leq1,\\[2pt]
 &u(x,y,0)=\frac{\partial u(x,y,0)}{\partial t}=0,~~(x,y)\in\bar{\Omega},\\[2pt]
 &u(x,y,t)=0,~~(x,y)\in\partial\Omega,~~0<t\leq 1.
 \end{array}$$
 Note that the equation can be equivalently written as
 $$\begin{array}{rl}
 \frac{\partial u(x,y,t)}{\partial t}={_{0}I^\alpha_t(\Delta u)}
 +\sin(x)\sin(y)\big[(\alpha+3)t^{\alpha+2}
 +\frac{2\Gamma(\alpha+4)}{\Gamma(2\alpha+4)}t^{2\alpha+3}\big],
 \end{array}$$
 where $\alpha=\gamma-1$. The exact solution for this problem is $u(x,t)=\sin(x)\sin(y)t^{\alpha+3}$.
 \end{example}

\begin{table}[hbt!]
\begin{center}
\caption{Numerical convergence orders in temporal direction with $h=\frac{\pi}{16}$ for Example \ref{ex1}.}\label{table1}
\renewcommand{\arraystretch}{0.9}
\def\temptablewidth{0.7\textwidth}
{\rule{\temptablewidth}{0.7pt}}
 \begin{tabular*}{\temptablewidth}{@{\extracolsep{\fill}}llcccc}
 $\alpha$ & $\tau$
 &\multicolumn{2}{c}{Proposed scheme (\ref{compact-scheme})}
 &\multicolumn{2}{c}{Compact scheme in \cite{Zhang2}}\\
 \cline{3-4}\cline{5-6}
 & & $E_{\infty}(h,\tau)$ &Rate1 &$E_{\infty}(h,\tau)$ &Rate1\\\hline
 $0.25$       & 1/5    & 6.9507e-3  & $\ast$  & 2.7048e-2   & $\ast$\\
              & 1/10   & 1.7717e-3  & 1.9720  & 8.4482e-3   & 1.6788\\
              & 1/20   & 4.4606e-4  & 1.9898  & 2.5877e-3   & 1.7070\\
              & 1/40   & 1.1292e-4  & 1.9819  & 7.8500e-4   & 1.7209\\
              & 1/80   & 2.8847e-5  & 1.9688  & 2.3742e-4   & 1.7253\\[4pt]
 $0.5$        & 1/5    & 1.0421e-2  & $\ast$  & 8.2035e-2   & $\ast$\\
              & 1/10   & 2.6014e-3  & 2.0021  & 2.9942e-2   & 1.4541\\
              & 1/20   & 6.5195e-4  & 1.9965  & 1.0749e-2   & 1.4780\\
              & 1/40   & 1.6294e-4  & 2.0004  & 3.8291e-3   & 1.4891\\
              & 1/80   & 4.1060e-5  & 1.9886  & 1.3594e-3   & 1.4941\\[4pt]
 $0.75$       & 1/5    & 1.7341e-2  & $\ast$  & 1.9340e-1   & $\ast$\\
              & 1/10   & 4.3653e-3  & 1.9901  & 8.1577e-2   & 1.2454\\
              & 1/20   & 1.0899e-3  & 2.0019  & 3.4379e-2   & 1.2466\\
              & 1/40   & 2.7235e-4  & 2.0007  & 1.4484e-2   & 1.2471\\
              & 1/80   & 6.8480e-5  & 1.9917  & 6.0986e-3   & 1.2479\\
 \end{tabular*}
 {\rule{\temptablewidth}{0.7pt}}
 \end{center}
 \end{table}

\begin{table}[hbt!]
\begin{center}
\caption{Numerical convergence orders in spatial direction with $\tau=\frac1{10000}$ when $\alpha=0.1$ for Example \ref{ex1}.}
 \label{table2}
\renewcommand{\arraystretch}{0.9}
\def\temptablewidth{0.6\textwidth}
{\rule{\temptablewidth}{0.7pt}}
 \begin{tabular*}{\temptablewidth}{@{\extracolsep{\fill}}lcccc}
 $h$ &\multicolumn{2}{c}{Proposed scheme (\ref{compact-scheme})}
 &\multicolumn{2}{c}{Compact scheme in \cite{Zhang2}}\\
 \cline{2-3}\cline{4-5}
 & $E_{\infty}(h,\tau)$ &Rate2 &$E_{\infty}(h,\tau)$ &Rate2\\\hline
 $\pi/4$  & 5.0651e-4    & $\ast$ & 5.0632e-4 & $\ast$\\
 $\pi/8$  & 3.1111e-5    & 4.0251 & 3.1104e-5 & 4.0249\\
 $\pi/16$ & 1.9371e-6    & 4.0054 & 1.9421e-6 & 4.0014\\
 $\pi/32$ & 1.2245e-7    & 3.9837 & 1.2814e-7 & 3.9218\\
 \end{tabular*}
 {\rule{\temptablewidth}{0.7pt}}
 \end{center}
 \end{table}

 \begin{figure}\label{figure1}
 \begin{center}
 \includegraphics[height=6.3cm,width=17.3cm]{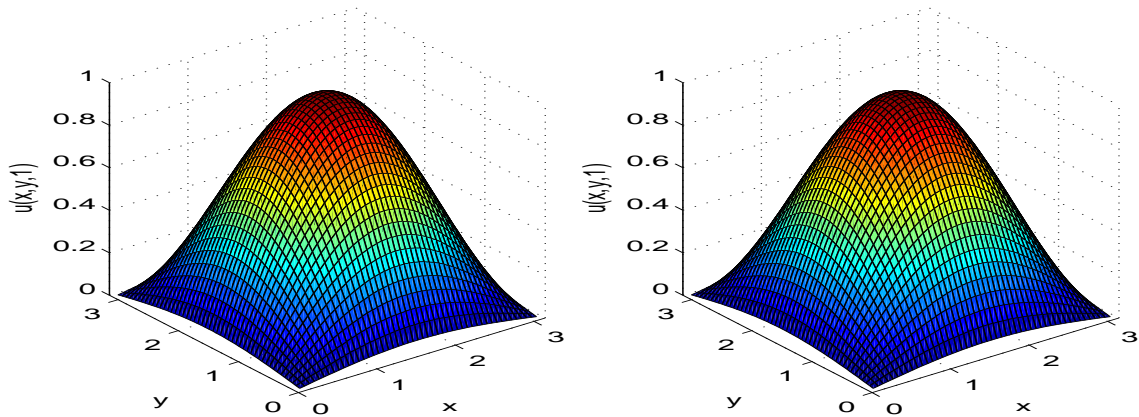}
 \end{center}
 \caption{The exact solution (left) and numerical solution (right) for Example \ref{ex1}, when $\alpha=0.5,~h=\tau=\frac1{50}$.}
 \end{figure}

 Our tests were done in MATLAB. At each time level,
 we need to solve a tridiagonal system with strictly diagonally dominant matrix.
 The structure of this matrix is same as that in \cite{Zhang2} and both schemes have similar computational cost.

 We now let $h_1=h_2=h$,
 the maximum norm errors between the exact and the numerical solutions
 $$E_\infty(h,\tau)=\max\limits_{0\leq k\leq N}
 \max\limits_{(x_i,y_j)\in\Omega}|u(x_i,y_j,t_k)-u_{ij}^k|$$
 are shown in Table \ref{table1} and  Table \ref{table2}. Furthermore, the temporal convergence order and spatial convergence order, denoted by
 $$Rate1=\log_2\bigg(\frac{E_\infty(h,2\tau)}{E_\infty(h,\tau)}\bigg)
 ~~\mbox{and}~~ Rate2=\log_2\bigg(\frac{E_\infty(2h,\tau)}{E_\infty(h,\tau)}\bigg),$$
 respectively, are reported. Table \ref{table1} shows that the proposed scheme (\ref{compact-scheme}) is more effective than the compact scheme in \cite{Zhang2}. These tables confirm the theoretical analysis. Meanwhile, Figure 1 shows the exact solution (left) and numerical solution (right) for Example \ref{ex1}, when $\alpha=0.5,~h=\tau=\frac1{50}$.

 \section*{Acknowledgment}
 The authors would like to express their gratitude to the referees for their valuable suggestions on the manuscript.

\end{document}